\documentclass[12pt]{article}
\usepackage[catalan,spanish,english]{babel}
\usepackage[psamsfonts]{amssymb}
\usepackage{cmmib57}
\usepackage{exscale}
\usepackage{amsmath}
\usepackage{amscd}
\usepackage{latexsym}
\usepackage{graphicx}
\newtheorem{thm}{Theorem}[section]

\newtheorem{lem}[thm]{Lemma}

\newcommand{\qed}{\quad{$\square$}}

\newcommand{\R}{\mathbb{R}}
\newcommand{\lga}{\longrightarrow}
\newcommand{\x}{\times}

\begin{document}
\begin{titlepage}
\null
\vspace{2cm}
\begin{center}
{\Large\bf Probability density for a hyperbolic SPDE  \\[2mm]
with time dependent coefficients }\\[2mm]
\medskip

by\\
\vspace{7mm}

\begin{tabular}{l@{\hspace{10mm}}l@{\hspace{10mm}}l}
{\sc Marta Sanz-Sol\'e}$\,^{(\ast)}$ &and &{\sc Iv\'an Torrecilla-Tarantino}$\,^{(\ast)}$\\
{\small marta.sanz@ub.edu }         &&{\small itorrecilla@ub.edu}\\
\end{tabular}
\begin{center}
{\small Facultat de Matem\`atiques}\\
{\small Universitat de Barcelona } \\
{\small Gran Via 585} \\
{\small 08007 Barcelona, Spain} \\
\end{center}
\end{center}

\vspace{2.5cm}

\noindent{\bf Abstract:} We prove the existence and smoothness of density for the
solution of a hyperbolic SPDE with free term coefficients depending on time, under hypoelliptic
non degeneracy conditions. The result extends those proved in \cite{CM02} to an infinite dimensional
setting.

\medskip

Two-parameter processes.

\vspace{2 cm}

\noindent

\footnotesize
{\begin{itemize}
\item[$^{(\ast)}$] Supported by the grant BMF 2003-01345 from the \textit{Direcci\'on
General de Investigaci\'on, Ministerio de Ciencia y Tecnolog\'{\i}a, Spain.}
\end{itemize}}
\end{titlepage}

\newpage

\section{Introduction}
\label{s1} The initial  developments of Malliavin Calculus
provided a probabilistic proof of H\"ormander's theorem for
hypoelliptic operators in square form. As an application,  the
existence and smoothness of the density for the solution of
diffusion processes with coefficients depending only on the
spatial variable were obtained (see \cite{Ma78}). There have been
several attempts to extend this result to diffusions with
coefficients depending on two variables, time and space. The first
results in this direction, \cite{Ta85}, \cite{CZ91}, apply to
pretty smooth coefficients (see the cases termed in Section 2 as
{\it smooth}, {\it factorable} and {\it regular H\"older}). More
recently, Cattiaux and Mesnager (\cite{CM02}) solved the problem
for hypoelliptic coefficients under less restrictive smoothness
conditions on the coefficients.

The  classical application of Malliavin Calculus  mentioned before
has been extended in \cite {NS85} to the two-parameter It\^o
equation -a wave equation in reduced form. Rules of two-parameter
It\^o calculus differ from those of the classical one. As a
consequence, the analogue of H\"ormander's condition is formulated
in terms of the {\it covariant derivative} instead of the Lie
brackets. In view of the results of \cite{CM02}, a natural
question is weather one could also extend the results from
\cite{NS85} to coefficients of the equation depending on the
two-dimensional time parameter. This article is devoted to study
this problem. Combining the techniques of \cite {NS85} with those
of \cite{CM02}, we  prove in Theorem \ref{t3.1} such an extension.

When studying the inverse of the Malliavin matrix corresponding to homogeneous diffusion processes,
one needs estimates of the type
$$
P\left\{\int_0^S Y_t^2 dt \le a \varepsilon^\delta, \int_0^S \alpha_t^2 dt\ge b \varepsilon^\eta\right\} \le \varepsilon ^p,
$$
for small $\varepsilon$ and bounded stopping time $S$. Here
$Y_t=Y_0+M_t+V_t$, $t\ge 0$, is a continuous semimartingale with
martingale and bounded variation components, $M_t$ and  $V_t$,
respectively; it is assumed that the quadratic variation of $M_t$
is of the form $\int_0^t \alpha_s^2 ds$ (see for instance
\cite{St83}, \cite{No86} for different proofs of this result).

When dealing with non-homogeneous diffusions with somehow rough coefficients,
such result does not suffice.
Actually, the crucial steps of the proofs of Theorem 2.6 and Proposition 3.2 in \cite{CM02} consists
of establishing alternate suitable extensions using a new approach an novel ideas.
Setting these ideas in an abstract framework, we
prove in Section \ref{s3} more sophisticated versions of Stroock-Norris type estimates
which are shown to be useful for two-parameter processes.

The paper is organized as follows. In Section \ref{s2}, we introduce the notation, the different hypothesis we are going to consider along the paper, and we state the main result. Section \ref{s3} is devoted to the extension of Stroock-Norris
estimates. With these tools, we prove in Section \ref{s4} the main result.

 \section{Notations, assumptions and the main result}
 \label{s2}

Consider the stochastic differential equation on $\mathbb{R}^m$
\begin{equation}
X_{z}=x_{0}+\int_{R_{z}} \left[\sum^{d}_{l=1}
A_{l}(r,X_r)\,dW^{l}_{r}+A_{0}(r,X_r)\,dr\right],\label{Itoeq}
\end{equation}
where $z=\left(s,t\right)\in[0,S]\times [0,T] $, $0\le
S,T<\infty$, $x_0\in\mathbb{R}^m$, $R_{z}=[0,s] \times [0,t]$ and
$W_{z}=\left(W^{1}_{z},\dots,W^{d}_{z}\right)$ is a
$d$-dimensional Brownian sheet (see \cite{CW75}).

Set  $\mathcal{E}=\left\{\left(s,t\right) \in R_{S,T}: \; st\neq 0\right\}$.
We assume that the coefficients of $(\ref{Itoeq})$ satisfy the
conditions
\begin{itemize}
\item[(h1)]$A_{l}:R_{S,T} \x \R^{m}\lga\R^{m}$, $0\leq l\leq d$,
are $\gamma$-H\"{o}lder continuous in $t$, for some $\gamma \in
\left]0,1\right[$, measurable  with respect to $s$, and infinitely
differentiable with respect to the spatial variable. Moreover,
\begin{equation}
K_{\gamma}:=\sup_{y\in \mathbb{R}^{m}}\sup_{0 \leq \theta \leq
S}\max_{0\leq l\leq
d}\left\|A_{l}\left(\theta,\cdot,y\right)\right\|_\gamma<\infty,
\label{ACOT2''}
\end{equation}
where the notation $\Vert \cdot\Vert_\gamma$ refers to
the usual H\"older norm.

\item[(h2)]for any multi-index
$\overline{\alpha}=\left(\alpha_{1},\ldots,\alpha_{m}\right)$, $\left|\overline{\alpha }\right|\geq 0$,
and  $0\leq l\leq d$, the
partial derivatives $\partial_{\overline{\alpha}}^{x}A_{l}$ with
respect to $x\in\mathbb{R}^{m}$ exist and
\begin{equation}
K:=\sup_{0\leq \theta \leq S}\sup_{0\leq \tau \leq T}\max_{0\leq
\left|\overline{\alpha }\right| \leq N+2}\max_{0\leq j\leq
d}\left\|
\partial_{\overline{\alpha}}^{x}A_{j}\left(\theta,\tau,\cdot\right)\right\|
_{\infty }<\infty,\label{ACOT2}
\end{equation}

\item[(h3)] for a fixed $z=(s,t)\in \mathcal{E}$, the vector fields $A_{l}$'s, $1\leq l\leq d$ satisfy the
restricted H\"{o}rmander's condition stated as follows:

The vector space spanned by the vector
fields $A_{1},\dots,A_{d}$, $A_{i}^{\nabla}A_{j}$, $1\leq
{i,j}\leq d$, $A_i^{\nabla}\left(A_j^{\nabla}A_k\right)$, $1\leq
{i,j,k}\leq d$,$\dots$,
$A_{i_{1}}^{\nabla}\left(\dots\left(A_{i_{n-1}}^{\nabla}A_{i_{n}}\right)\dots\right)$,
$1\leq i_1,\dots,i_n\leq d$,$\dots$, at the point
$(0,t,x_0)$ has full rank.

Here, the notation $A_i^\nabla A_j$ denotes the covariant derivative of the vector
field $A_j$ along $A_i$.
\end{itemize}

Assumption (h2) implies that the coefficients of the equation are
Lipschitz functions and have linear growth in the spatial
variable, uniformly in the time variables. By the usual method of
Picard's iterates, one can prove existence and uniqueness of
solution for the equation (\ref{Itoeq}) and moreover, the solution
has almost surely continuous paths (see Lemma 3.1 in \cite{NS85}).
\medskip

Condition (h3) clearly implies the following. There exist $N:=N(x_{0},t)\in \mathbb{N}$, and positive real numbers $c_N:=c_{N}(x_{0},t)$,
$s_0:=s_{0}(t)$, $R:=R(t)$ such that for any $v\in S^{m-1}$,
\begin{equation}
\sum_{k=0}^{N}\sum_{V\in { \Sigma }_{k}}\langle
v,V(\theta,t,y)\rangle^{2}\geq c_N,\label{Hör''}
\end{equation}
for any $(\theta,y)\in [0,s_{0}] \x B(x_{0},R)$, where
$\Sigma _{0}=\{A_{l},\;1\leq l\leq d\}$,
$\Sigma_{k+1}=\{ A_{l}^{\nabla}V,\;1\leq l\leq d,\;V\in
\Sigma _{k}\}$.

We consider different  combinations of regularity of the coefficients and non degeneracy conditions of the underlying differential operator,
as follows.

\begin{enumerate}
\item \textbf{Elliptic case:} (h1) and (h2) holds. Moreover,
the vector space spanned by the vector
fields $A_{1},\dots,A_{d}$  at the point $(0,t,x_0)$ has full rank.


\item \textbf{Smooth case:} (h1) to (h3) hold. In addition, for each $0\leq l\leq d$, the functions  $\partial_{\overline{\alpha
}}^{x}A_{l}$'s are $C_b^1$ in $s$,  for all multi-index
$\overline{\alpha }$,
\begin{equation}
K_{1}:=\sup_{0\leq \theta \leq S}\sup_{0\leq \tau \leq
T}\max_{0\leq \left|\overline{\alpha }\right| \leq N}\max_{1\leq
j\leq
d}\left\|\partial^{\theta}\partial_{\overline{\alpha}}^{x}A_{j}\left(\theta,\tau,\cdot\right)\right\|
_{\infty } < \infty. \label{ACOT22}
\end{equation}

\item \textbf{Factorable case:}
$A_{l}\left(\theta,\tau,x\right)=f_{l}\left(\theta\right)\overline{A}_{l}\left(\tau,x\right)$,
$0\leq l\leq d$, with $f_{l}$ measurable,
 $\frac{1}{c'} \geq |f_{l}| \geq c'>0$.  The functions  $\overline{A}_{l}$'s,
$0\leq l\leq d$, satisfy the analogue of hypothesis (h1)-(h3) for
coefficients which do not depend on $\theta$.

\item \textbf{Regular H\"{o}lder case:} (h1) to (h3) hold. In addition, for any multi-index $\overline{\alpha }$, the functions
$\partial_{\overline{\alpha }}^{x}A_{j}$, $0\leq j\leq d$,  are
$\beta\left(\overline{\alpha }\right)$-H\"{o}lder continuous in the argument
$s$, for some $\beta\left(\overline{\alpha }\right)\in\left]\frac{1}{2},1\right[$. Moreover,
\begin{equation}
K_{\beta(\overline\alpha)}:=\sup_{y \in \mathbb{R}^{m}}\sup_{0 \leq \tau \leq
T}\max_{0\leq \left|\overline{\alpha}\right| \leq N+2}\max_{0\leq l\leq d}
\left\|\partial_{\overline{\alpha}}^{x}A_{l}\left(\cdot,\tau,y
\right)\right\|_{\beta\left(\overline{\alpha}\right)}<\infty.\label{ACOT1}
\end{equation}

\item \textbf{Irregular H\"{o}lder case:} The same assumptions as in the preceding case, except that here
$\beta\left(\overline{\alpha }\right)\in \left]0,\frac{1}{2}\right]$.
\end{enumerate}


The main result of this paper is the next theorem, stating
the existence and smoothness of density for the probability
law of the solution of (\ref{Itoeq}) at any fixed point $z\in \mathcal{E}$.

\begin{thm}
\label{t3.1}
Let $X=\left\{X_{z},\; z \in R_{S,T}\right\}$ be the solution of (\ref{Itoeq}).
Each one of the  set of assumptions termed before as elliptic, smooth, factorable, regular H\"older and irregular H\"older, imply that  the random vector
$X_{z}$, for fixed $z\in\mathcal{E}$, has an infinitely differentiable density with respect to the Lebesgue measure.
\end{thm}
{\bf Remark} In the formulation of assumptions (h1)-(h3) and of the different scenaries, the roles of the time components $s$ and $t$ might be exchanged.

\section{Stroock-Norris type lemmas for continuous semimartingales depending on a parameter}
\label{s3}

In this section, we prove two extensions of the Stroock-Norris
estimates. The difference between them stands on the order of
H\"older continuity of the involved processes. We notice that
Lemma \ref{l3.1} could provide as a by-product an alternate proof
of Norris Lemma (\cite{No86})).

\begin{lem}
\label{l3.1}
Let $(Y_s(\lambda), s\in[0,S])$ be a real continuous
semimartingale depending on a parameter $\lambda\ge 0$,  with decomposition
\begin{equation}
Y_s\left(\lambda\right)=Y_{0}\left(\lambda\right)+M_s\left(\lambda\right)+V_s\left(\lambda\right),
\label{martrep2}
\end{equation}
where 
$M_s\left(\lambda\right)$, $V_s\left(\lambda\right)$ denote a
local martingale and a bounded variation process, respectively, satisfying
\begin{align*}
M_s\left(\lambda\right)&=\sum_{j=1}^{m}\int_{0}^s\Psi_{\eta}^{j}
\left(\lambda\right)d\widetilde{M}_{\eta}^{j},\quad
\left\langle \widetilde{M}^{j},\widetilde{M}^{k}
\right\rangle_s
=\int_{0}^s\Theta_{\eta}^{j,k}d\eta, \\
V_s\left(\lambda\right)&=\int_{0}^s\Phi_{\eta}\left(\lambda\right)d\eta.
\end{align*}
Assume that:

(i) For each $\lambda\ge 0$, $\Psi_{\eta}^{j}\left(\lambda\right)$,
$\Phi_{\eta}\left(\lambda\right)$ and $\Theta_{\eta}^{j,k}$, $1\leq j,k \leq m$,  are adapted  continuous
processes, indexed by $\eta\in[0,S]$, bounded by some constant $\mathcal{K}$,
uniformly in $\eta,\lambda$.

(ii) For each $\eta\in [0,S]$, $1\leq j \leq m$, $Y_{\eta}(\lambda)$,
$Y_{0}\left(\lambda\right)$, $\Psi_{\eta}^{j}\left(\lambda\right)$,
 $\Phi_{\eta}\left(\lambda\right)$ as functions of $\lambda$, are
$\beta$-H\"{o}lder continuous, with $\frac{1}{2}<\beta<1$, uniformly in $\eta$.


Set
$
\left \langle M\left(\lambda\right) \right
\rangle_s=\int_{0}^s\Upsilon_{\eta}\left(\lambda\right)d\eta,
$
where
\begin{equation}
\Upsilon_{\eta}\left(\lambda\right)=
\sum_{j,k=1}^m \Psi_{\eta}^{j}\left(\lambda\right)
\Psi_{\eta}^{k}\left(\lambda\right)\Theta_{\eta}^{j,k}.
\label{defPsi}
\end{equation}
Fix $\nu>\frac{3}{2\beta-1}$.
Then, for any
$\rho>3+2\nu$, positive constants $\alpha_{1}$, $\alpha_{2}$, $p \geq 2$,
and $\varepsilon$ small enough, there exists a constant $C$ such that
\begin{equation*}
\mathbf{P}\left\{\int_{0}^{s}Y_{u}^{2}\left(u\right)du \leq
\alpha_{1}\varepsilon^{\rho},\;\int_{0}^{s}\Upsilon_{u}\left(u\right)du
\geq \alpha_{2}\varepsilon\right\} \leq C \varepsilon^{p}.
\end{equation*}
\end{lem}

{\it Proof}: 
Fix $n\ge 1$ and set $s_{i}=\frac{is}{n}$, $i=0,\ldots,n$,
For each $i=0,\ldots ,n-2$, we define
\begin{equation*}
D_{i}=\left\{\int_{0}^{s}Y_{u}^{2}\left(u\right)du\leq
\alpha_{1}\varepsilon^{\rho},\;\int_{s_{i}}^{s_{i+1}}\Upsilon_{u}\left(u\right)du
\geq \frac{\alpha_{2}\varepsilon}{2\left(n-1\right)}\right\}.
\end{equation*}
We also define
\begin{equation*}
D_{n-1}=\left\{\int_{0}^{s}Y_{u}^{2}\left(u\right)du\leq
\alpha_{1}\varepsilon^{\rho},\;\int_{\left(1-\frac{1}{n}\right)s}^{s}\Upsilon_{u}\left(u\right)du
\geq \frac{\alpha_{2}\varepsilon}{2}\right\}.
\end{equation*}
Then
\begin{equation}
\left\{\int_{0}^{s}Y_{u}^{2}\left(u\right)du \leq
\alpha_{1}\varepsilon^{\rho},\;\int_{0}^{s}\Upsilon_{u}\left(u\right)du
\geq \alpha_{2}\varepsilon\right\}\subset
\bigcup_{i=0}^{n-1}D_{i}. \label{Inclusiondiscret}
\end{equation}
Chebyshev's inequality and the boundedness of $\Psi^{j}$,
$\Theta^{j,k}$, $1\leq j,k \leq m$, yield for all $p' \geq 1$,
\begin{equation*}
\mathbf{P}\left\{D_{n-1}\right\} \leq
\mathbf{P}\left\{\int_{\left(1-\frac{1}{n}\right)s}^{s}\Upsilon_{u}\left(u\right)du
\geq \frac{\alpha_{2}\varepsilon}{2}\right\}\leq
\left(\frac{2m^{2}\mathcal{K}^{3}s}{\alpha_{2}}\right)^{p'}\frac{\varepsilon^{-p'}}{n^{p'}}.
\end{equation*}
Taking $n=\left[ {\varepsilon}^{-\nu}\right]$ (where $\left[ \cdot
\right]$ denotes the integer part) and $\nu > 1$, we have
$\frac{\varepsilon^{-p'}}{n^{p'}}\le
2^{p'}\varepsilon^{\left(\nu-1\right)p'}$, for any $\varepsilon <
\varepsilon_{0}=2^{-\frac{1}{\nu}}$. Therefore,
\begin{equation}
\mathbf{P}\left\{D_{n-1}\right\}\leq
C\varepsilon^{p''},\label{LastSetEstimate}
\end{equation}
for all $p''=\left(\nu-1\right)p' \geq 2$.

We now study the terms $\mathbf{P}(D_{i})$, for all $i=0,\ldots,n-2$.
Setting
\begin{equation*}
F_{i}=\left\{\int_{s_{i}}^{s_{i+2}}Y_{u}^{2}\left(u\right)du\leq
\alpha_{1}\varepsilon^{\rho},\;\int_{s_{i}}^{s_{i+1}}\Upsilon_{u}\left(u\right)du
\geq \frac{\alpha_{2}\varepsilon}{2n}\right\},
\end{equation*}
we clearly have $D_{i}\subset F_{i}$. 

By It\^{o}'s formula, for $u\in
\left]s_{i},s_{i+2}\right]$, it holds that
\begin{align}
Y_{u}^{2}\left(\lambda\right)&= Y_{s_{i}}^{2}\left(\lambda\right)
+2\sum_{j=1}^{m}\int_{s_{i}}^{u}Y_{\eta}\left(\lambda\right)\Psi_{\eta}^{j}
\left(\lambda\right)d\widetilde{M}_{\eta}^{j}\nonumber\\
&+2\int_{s_{i}}^{u}Y_{\eta}\left(\lambda\right)\Phi_{\eta}\left(\lambda\right)d\eta+
\int_{s_{i}}^{u}\Upsilon_{\eta}\left(\lambda\right)d\eta.\label{ItoHolder}
\end{align}

Set $\lambda=u$ in $(\ref{ItoHolder})$. By integrating with respect to the
$u$ variable and applying
Fubini's theorem and a stochastic Fubini theorem, we obtain
\begin{equation}
\int_{s_{i}}^{s_{i+2}}Y_{u}^{2}\left(u\right)du \geq
\mathcal{C}_{s_{i+2}}+\mathcal{A}_{s_{i+2}}+
\mathcal{M}_{s_{i+2}},  \label{estimacioHC}
\end{equation}
where, for any $s\ge s_i$,
\begin{eqnarray*}
\mathcal{C}_s&=&
\int_{s_{i}}^{s}\left(\int_{\eta}^{s}\Upsilon_{\eta}\left(u\right)du\right)
d\eta, \\
\mathcal{A}_{s}&=&2\int_{s_{i}}^{s}\left(\int_{\eta}^{s}
Y_{\eta}\left(u\right)\Phi_{\eta}\left(u\right)du\right)
d\eta, \\
\mathcal{M}_{s}&=&2\sum_{j=1}^{m}\int_{s_{i}}^{s}\left(\int_{\eta}^{s}
Y_{\eta}\left(u\right)\Psi_{\eta}^{j} \left(u\right)du\right)
d\widetilde{M}_{\eta}^{j}.
\end{eqnarray*}
We  devote the remaining of the proof to show the inclusion
\begin{align}
&F_{i}\subset \left\{\sup_{s_{i} \leq u \leq s_{i+2}}\left|
\mathcal{M}_{u}\right| \geq
\frac{\alpha_{2}s}{16}\varepsilon^{1+2\nu},\;\left\langle\mathcal{M}\right\rangle_{s_{i+2}}
< \frac{\alpha_{2}^{2}s^{2}}{256}{\varepsilon}^{ 3+4\nu }\right\}.
\label{maininclusionHC}
\end{align}
for some $\nu>0$.  Then, applying the martingale exponential
inequality,
\begin{equation*}
\mathbf{P}\left\{\sup_{s_{i} \leq u \leq s_{i+2}}\left|
\mathcal{M}_{u}\right| \geq
\frac{\alpha_{2}s}{16}\varepsilon^{1+2\nu},\;\left\langle\mathcal{M}\right\rangle_{s_{i+2}}
< \frac{\alpha_{2}^{2}s^{2}}{256}{\varepsilon}^{ 3+4\nu }\right\}
\leq 2\exp \left(-\frac{1}{2}\varepsilon^{-1}\right),
\end{equation*}
and this will finish the proof.

The above inclusion is obtained by proving a lower bound for
$\mathcal{C}_{s_{i+2}}$, and upper bounds for $\left|
\mathcal{A}_{s_{i+2}}\right|$ and $\left\langle
\mathcal{M}_{\cdot}\right\rangle_{s_{i+2}}$.
\smallskip

{\emph{Lower bound for }$\mathcal{C}_{s_{i+2}}$.} Clearly, on
$F_{i}$,  the triangular inequality and the H\"older continuity of
$\Upsilon_{\eta}(\cdot)$ (with constant $\mathcal{K}_{\beta}$), imply
\begin{align*}
\mathcal{C}_{s_{i+2}}
 &\geq\int_{s_{i}}^{s_{i+2}}\left(\int_{\eta}^{s_{i+2}}\Upsilon_{\eta}\left(\eta\right)
du\right)d\eta-
\int_{s_{i}}^{s_{i+2}}\left(\int_{\eta}^{s_{i+2}}\left|\Upsilon_{\eta}\left(u\right)-
\Upsilon_{\eta}\left(\eta\right)\right|du\right)d\eta\\
&\geq\int_{s_{i}}^{s_{i+1}}\left(s_{i+2}-\eta\right)\Upsilon_{\eta}\left(\eta\right)
d\eta-\int_{s_{i}}^{s_{i+2}}\left(\int_{\eta}^{s_{i+2}}\mathcal{K}_{\beta}\left|u-\eta\right|^{\beta}
du\right)d\eta\\
&\geq\left(s_{i+2}-s_{i+1}\right)\int_{s_{i}}^{s_{i+2}}\Upsilon_{\eta}\left(\eta\right)
d\eta-\mathcal{K}_{\beta}\left(s_{i+2}-s_{i}\right)^{\beta+2}\\
&\geq \frac{\alpha_{2}s\varepsilon}{2n^{2}}-
\mathcal{K}_{\beta}\frac{\left(2s\right)^{\beta+2}}{n^{\beta+2}}.
\end{align*}
Taking $n=\left[ {\varepsilon}^{-\nu}\right]$  and $\nu > 0$, yields
\begin{equation*}
\frac{\alpha_{2}s\varepsilon}{2n^{2}}-
\mathcal{K}_{\beta}\frac{\left(2s\right)^{\beta+2}}{n^{\beta+2}}\ge
\frac{\alpha_{2}s}{2}\varepsilon^{1+2\nu}- \mathcal{K}_{\beta}(4s
\varepsilon^\nu)^{\beta+2},
\end{equation*}
for any $\varepsilon < \varepsilon_0$. Choose $\nu>0$ such that
$\nu\beta>1$, so that $\nu\left(\beta+2\right)>1+2\nu$. Thus,
taking $ \varepsilon_{1}=\left(
\frac{\alpha_{2}}{\mathcal{K}_{\beta}
4^{\beta+3}s^{\beta+1}}\right)^{\frac{1}{\nu\beta-1}}, $ we have
that for all $\varepsilon<\varepsilon_{1}\wedge \varepsilon_0$, on
$F_i$,
\begin{equation}
\mathcal{C}_{s_{i+2}}\geq \frac{\alpha_{2}s}{4}{\varepsilon}^{
1+2\nu }. \label{CepsilonHC}
\end{equation}

{\emph{Upper bound for }$\left|
\mathcal{A}_{s_{i+2}}\right| $.}  Jensen's inequality implies
\begin{equation*}
\left|\mathcal{A}_{s_{i+2}}\right|\leq 2\mathcal{K}
\int_{s_{i}}^{s_{i+2}}\left(\int_{\eta}^{s_{i+2}}
\left|Y_{\eta}\left(u\right)\right|du\right) d\eta.
\end{equation*}
Since $\lambda \to Y_{u}\left(\lambda\right)$ is $\beta$-H\"{o}lder
continuous (H\"older constant $\mathcal{K}'_\beta$), we obtain
\begin{align*}
&\int_{s_{i}}^{s_{i+2}}\left(\int_{\eta}^{s_{i+2}}
\left|Y_{\eta}\left(u\right)\right|du\right) d\eta\\
& \leq \int_{s_{i}}^{s_{i+2}}\left(\int_{\eta}^{s_{i+2}}
\left|Y_{\eta}\left(\eta\right)\right|du\right)
d\eta+\int_{s_{i}}^{s_{i+2}}\left(\int_{\eta}^{s_{i+2}}
\left|Y_{\eta}\left(u\right)-Y_{\eta}\left(\eta\right)\right|du\right)
d\eta\\
&\leq \int_{s_{i}}^{s_{i+2}}\left(s_{i+2}-\eta \right)
\left|Y_{\eta}\left(\eta\right)\right|
d\eta+\mathcal{K}'_{\beta}\int_{s_{i}}^{s_{i+2}}\left(\int_{\eta}^{s_{i+2}}
\left|u-\eta\right|^{\beta}du\right) d\eta\\
&\leq \int_{s_{i}}^{s_{i+2}}\left(s_{i+2}-\eta \right)
\left|Y_{\eta}\left(\eta\right)\right|
d\eta+\mathcal{K}'_{\beta}\left(s_{i+2}-s_{i}\right)^{\beta+2}.
\end{align*}
But,
\begin{align*}
&\int_{s_{i}}^{s_{i+2}}\left(s_{i+2}-\eta \right)
\left|Y_{\eta}\left(\eta\right)\right| d\eta
\leq\left\{\int_{s_{i}}^{s_{i+2}}\left(s_{i+2}-\eta
\right)^{2}d\eta\right\}^{\frac{1}{2}}\left\{\int_{s_{i}}^{s_{i+2}}
Y_{\eta}^{2}\left(\eta\right)d\eta\right\}^{\frac{1}{2}}\\
&\leq\left(s_{i+2}-s_{i}\right)^{\frac{3}{2}}\left\{\int_{s_{i}}^{s_{i+2}}
Y_{\eta}^{2}\left(\eta\right)d\eta\right\}^{\frac{1}{2}}.
\end{align*}
Thus, on $F_{i}$, we have
\begin{align*}
\left|\mathcal{A}_{s_{i+2}}\right|\leq
2\mathcal{K}\left(\left(\frac{2s}{n}\right)^{\frac{3}{2}}{\alpha_{1}}^\frac{1}{2}\varepsilon^{\frac{\rho}{2}}
+\mathcal{K}^\prime_\beta
\left(\frac{2s}{n}\right)^{\beta+2}\right).
\end{align*}
As before, taking $n=\left[ {\varepsilon}^{-\nu}\right]$, $\nu>0$,
one obtains
\begin{align*}
\left(\frac{2s}{n}\right)^{\frac{3}{2}}{\alpha_{1}}^\frac{1}{2}\varepsilon^{\frac{\rho}{2}}
+\mathcal{K}^\prime_\beta \left(\frac{2s}{n}\right)^{\beta+2}&\le
(4s)^{\frac{3}{2}}{\alpha_{1}}^\frac{1}{2}\varepsilon^{\frac{3\nu}{2}+\frac{\rho}{2}}
+\mathcal{K}^\prime_\beta\left(4s\right)^{\beta+2}\varepsilon^{\nu\left(\beta+2\right)}\\
&\leq \left(\left(4s\right)^{\frac{3}{2}}{\alpha_{1}}^\frac{1}{2}
+\mathcal{K}^\prime_\beta\left(4s\right)^{\beta+2}\right)
{\varepsilon}^{(\frac{3\nu}{2}+\frac{\rho}{2})\wedge
(\nu(\beta+2))},
\end{align*}
for any $\varepsilon < \varepsilon_0$. Choose $\nu$ such that $\nu
\beta>1$ and $\rho>2+\nu$. Then $
m\left(\rho,\nu,\beta\right)=\left\{\frac{\rho}{2}+\frac{3\nu}{2}\right\}\wedge
\left\{\nu\left(2+\beta\right)\right\}>1+2\nu. $ Thus, taking
\begin{equation*}
\varepsilon^{\prime}_{1}=\left(
\frac{\alpha_2}{32\mathcal{K}\left(4^{\frac{3}{2}}s^{\frac{1}{2}}{\alpha_{1}}^\frac{1}{2}
+\mathcal{K}^\prime_\beta
4^{\beta+2}s^{\beta+1}\right)}\right)^{\frac{1}{m\left(\rho,\nu,\beta\right)-1-2\nu}},
\end{equation*}
on $F_i$, we obtain for all $\varepsilon
<\varepsilon^{\prime}_{1}\wedge\varepsilon_{0}$,
\begin{equation}
\left| \mathcal{A}_{s_{i+2}}\right| \leq
\frac{\alpha_{2}s}{16}{\varepsilon}^{ 1+2\nu }.
\label{AepsilonHC}
\end{equation}

{\emph{Upper bound for }$\left\langle
\mathcal{M}\right\rangle_{s_{i+2}}$}. Clearly,
\begin{align*}
\left\langle\mathcal{M}\right\rangle_{s_{i+2}}&=4\sum_{j,k=1}^{m}\int_{s_{i}}^{s_{i+2}}\left(\int_{\eta}^{s_{i+2}}
Y_{\eta}\left(u\right)\Psi_{\eta}^{j}
\left(u\right)du\right)\\
&\quad\times\left(\int_{\eta}^{s_{i+2}}
Y_{\eta}\left(u\right)\Psi_{\eta}^{k} \left(u\right)du\right)
\Theta_{\eta}^{j,k}d\eta.
\end{align*}
Jensen's inequality implies
\begin{align*}
\left|\int_{\eta}^{s_{i+2}} Y_{\eta}\left(u\right)\Psi_{\eta}^{j}
\left(u\right)du\right|\leq \mathcal{K}
\left(s_{i+2}-\eta\right)\left|Y_{\eta}\left(\eta\right)\right|
+\mathcal{K}\mathcal{K}'_{\beta}\left(s_{i+2}-\eta\right)^{\beta+1}.
\end{align*}
Thus,
\begin{equation*}
\left\langle\mathcal{M}\right\rangle_{s_{i+2}}
\leq
C_1\left(s_{i+2}-s_{i}\right)^{2}\int_{s_{i}}^{s_{i+2}}
Y_{\eta}^{2}\left(\eta\right)d\eta
+C_2
\left(s_{i+2}-s_{i}\right)^{2\beta+3},
\end{equation*}
and on $F_{i}$,
\begin{equation*}
\left\langle\mathcal{M}\right\rangle_{s_{i+2}}\leq
C_1\left(\frac{2s}{n}\right)^{2} \alpha_{1}\varepsilon^{\rho}
+C_2
\left(\frac{2s}{n}\right)^{2\beta+3}.
\end{equation*}

Proceeding as before, we may choose $\nu>0$ such that
$\nu(2\beta-1)>3$, then $\rho>3+2\nu$, and find
$\varepsilon_1^{\prime\prime}$ such that on $F_i$, for all
$\varepsilon <\varepsilon^{\prime \prime}_{1}$,
\begin{equation}
\left\langle\mathcal{M}\right\rangle_{s_{i+2}} <
\frac{\alpha_{2}^{2}s^{2}}{256}{\varepsilon}^{ 3+4\nu}.
\label{QVMepsilonHC}
\end{equation}


The set
$
F_{i}\cap \left\{ \sup_{s_{i} \leq u \leq s_{i+2}}\left|
\mathcal{M}_{u}\right|
<\frac{\alpha_{2}s}{16}\varepsilon^{1+2\nu}\right\}
$
is empty. In fact, on this set, by
$\left(\ref{estimacioHC}\right)$, $\left(\ref{CepsilonHC}\right)$,
$\left(\ref{AepsilonHC}\right)$ and  (\ref{QVMepsilonHC}), we obtain
\begin{equation}
\label{11}
\alpha_{1}\varepsilon^{\rho} \geq
\int_{s_{i}}^{s_{i+2}}Y^{2}_{u}\left(u\right)du \geq
\mathcal{C}_{s_{i+2}}-\left|\mathcal{A}_{s_{i+2}}\right|-
\left|\mathcal{M}_{s_{i+2}}\right|\geq
\frac{\alpha_{2}s}{8}\varepsilon^{1+2\nu}.
\end{equation}
Hence, taking
$\beta >\frac{1}{2}$,
$\nu >\frac{3}{2\beta-1}$,
$\rho>(3+2\nu)$,
$\varepsilon <\tilde\varepsilon_0$, with
$$\tilde\varepsilon_0:=\min \left\{ \left(
\frac{\alpha_{2}s}{8\alpha_{1}}\right)^\frac{1}{\rho-1-2\nu},\varepsilon_0,\varepsilon
_{1},\varepsilon'_{1},\varepsilon''_{1}\right\},$$
(\ref{11}) cannot be satisfied.

Consequently, for any $\varepsilon < \tilde\varepsilon_0$,
\begin{align*}
F_i=&F_{i}\cap \left\{ \sup_{s_{i} \leq u \leq s_{i+2}}\left|
\mathcal{M}_{u}\right| \geq
\frac{\alpha_{2}s}{16}\varepsilon^{1+2\nu}\right\}\\
&\subset \left\{\sup_{s_{i} \leq u \leq s_{i+2}}\left|
\mathcal{M}_{u}\right| \geq
\frac{\alpha_{2}s}{16}\varepsilon^{1+2\nu},\;\left\langle\mathcal{M}\right\rangle_{s_{i+2}}
< \frac{\alpha_{2}^{2}s^{2}}{256}{\varepsilon}^{ 3+4\nu }\right\}.
\end{align*}
Thus, we obtain $\left(\ref{maininclusionHC}\right)$, and this ends the proof of the lemma.
\hfill\qed
\medskip

The next lemma applies to the more irregular situation where $\beta\in]0,\frac{1}{2}]$.
\begin{lem}
\label{l3.2}
 Consider a continuous semimartingale $(Y_s(\lambda), 0\le s\le S)$, with decomposition given in
 (\ref{martrep2}) and satisfying the assumption (i) of  Lemma \ref{l3.1}.
 Fix $\beta\in]0,\frac{1}{2}]$ and assume that:

(ii') For each $\eta\in[0,S]$, $1\leq j \leq m$,
$Y_{0}(\lambda)$,
$\Psi_{\eta}^{j}(\lambda)$,
$\Phi_{\eta}(\lambda)$, as functions of $\lambda$, are $\beta$-H\"{o}lder continuous in $\lambda$.

(iii') for all $\beta'<\beta $,
there exists versions of
$\Psi^{j}_{\eta}(\eta)$ and $\Theta^{j,k}_{\eta}$,
$1 \leq j,k \leq m$, respectively, which are $\beta'$-H\"{o}lder
continuous on $\left[0,s\right]$.

 Furthermore, for all $p\geq 2$
 \begin{equation*}
 \mathbb{E}\left(\left\|\Psi^{j}_{\cdot}\left(\cdot\right)\right\|^p_{\beta'}+
 \left\|\Theta^{j,k}_{\cdot}\left(\cdot\right)\right\|^p_{\beta'}\right)\le C_p<\infty.
 \end{equation*}


Then, for any
$\rho>\left(\frac{11}{2}+\frac{4}{\beta'}\right)\left(1+\frac{1}{\beta'}\right)$, positive constants
$\alpha_{1}$, $\alpha_{2}$, $p\geq 2$ and
 $\varepsilon$ sufficiently small, there exists a constant $C$ such that
\begin{equation*}
\mathbf{P}\left\{\int_{0}^{s}Y_{u}^{2}\left(u\right)du \leq
\alpha_{1}\varepsilon^{\rho},\;\int_{0}^{s}\Upsilon_{u}\left(u\right)du
\geq \alpha_{2}\varepsilon\right\} \leq C \varepsilon^{p},
\end{equation*}
\end{lem}

Proof: We follow the arguments of Proposition $3.2$ in
\cite{CM02} in a more abstract presentation.

{\bf Fact 1.}
{\sl
Fix $\beta'<\beta$,
$\varepsilon <s \wedge 1$,
$a>2\left(1+\frac{1}{\beta ^{\prime }}\right)$, $p\ge 2$.
Condition (ii') implies the existence of a positive constant
$c_{p}$ such that
\begin{equation*}
\mathbf{P}\left\{ \int_{0}^{s}Y_{u}^{2}\left(u\right)du \leq
\alpha_{1}\varepsilon^{a},\;\sup_{0\leq u \leq s}\left|
Y_{u}\left(u\right) \right| > \left(1+
\sqrt{\alpha_{1}}\right)\varepsilon \right\} \leq c_{p}\varepsilon
^{p},
\end{equation*}
}
This can be checked following the proof of
 Lemma 3.4 on \cite{CM02} with $\langle \xi,\Phi_u^{*-1} Y\rangle(x):=Y_u\left(u\right)$.

As a consequence, taking $\rho=ra$,
$\varepsilon^{r}<s \wedge 1$, with $r>0$ and $a$ as before,
we reduce the problem to estimate the probability of the set

\begin{equation*}
B=\left\{\sup_{0\leq u \leq s}\left| Y_{u}\left(u\right) \right|
\leq \left(1+
\sqrt{\alpha_{1}}\right)\varepsilon^{r},\;\int_{0}^{s}\Upsilon_{u}\left(u\right)du
\geq \alpha_{2}\varepsilon \right\}.
\end{equation*}

As in the previous Lemma \ref{l3.1}, we consider the
discretization of the interval $\left[0,s \right]$ given by
$s_{i}=\frac{is}{n}$, $i=0,\ldots ,n$,
with $n=\left[ \varepsilon ^{-\frac{4}{\beta ^{\prime }}}\right]$.

For $i=0,\dots,n-1$, set
\begin{eqnarray*}
C_{i}&:=&\left\{\exists u \in \left[s_{i},s_{i+1}\right]:\;
\Upsilon_{u}\left(u\right)< \alpha_{2}\varepsilon^{\frac{3}{2}},\;
\int_{0}^{s}\Upsilon_{u}\left(u\right)du \geq
\alpha_{2}\varepsilon
\right \},\\
D_{i}&:=&\left\{ \sup_{s_{i}\leq u \leq s_{i+1}}\left|
Y_{u}\left(u\right) \right| \leq \left(1+
\sqrt{\alpha_{1}}\right)\varepsilon^{r},\;
\Upsilon_{u}\left(u\right) \geq
\alpha_{2}\varepsilon^{\frac{3}{2}},\;\forall u \in
\left[s_{i},s_{i+1}\right]\right\}.
\end{eqnarray*}
Clearly, $B\cap D_{i}^{c} \subset C_{i}$.
\smallskip

{\bf Fact 2.} {\sl Condition (iii') implies that for any $p\ge 2$,
and $\varepsilon^{\frac{1}{2}}<\frac{1}{2s}$,
\begin{equation*}
\mathbf{P}\left( \bigcap_{i=0}^{n-1}C_{i}\right) \leq
c_{p}\varepsilon ^{p}.
\end{equation*} }

For the proof of this fact, we can follow the arguments of Lemma 3.7 in \cite{CM02}
with $\sum_{j=1}^m \langle \xi,\Phi_u^{*-1}[X_j,Z]\rangle^2 (x):=\Upsilon_u(u)$.

Let $D=\bigcup_{i=0}^{n-1}D_{i}$,
According to Fact 2, we have
$ \mathbf{P}\left( B\cap D^{c}\right) \leq c_{p}\varepsilon ^{p}$
for any $\varepsilon ^{\frac{1}{2}}<\frac{1}{2s}$.

Therefore, we have only to estimate $\mathbf{P}\left( B\cap D
\right) $. Moreover, since
\begin{equation*}
\mathbf{P}\left( B\cap D \right)\leq \mathbf{P}\left( D \right)
\leq n\;\underset{i}{\max}\left\{\mathbf{P}\left( D_{i}
\right)\right\} \leq \varepsilon ^{-\frac{4}{\beta'}}
\underset{i}{\max}\left\{\mathbf{P}\left( D_{i} \right)\right\},
\end{equation*}
it suffices  to show that
$
\mathbf{P}\left( D_{i} \right) \leq C_{p}\varepsilon ^{p},
$
for any $p\ge 1$, with some constant $C_{p}$, not depending on $i$.

Recall the definition of $\Upsilon$ given by
$\left(\ref{defPsi}\right)$. Owning to assumptions (i) and (ii'),
\begin{equation}
\sup_{\left(u,\lambda \right)\in \left[
s_{i},s_{i+1}\right]^{2}}\left|\Upsilon_{u}\left(u\right)
-\Upsilon_{u}\left(\lambda\right) \right| \leq C
\left|s_{i+1}-s_{i}\right|^{\beta}=C
\left(\frac{s}{n}\right)^\beta \le C
s^{\beta}\varepsilon^{\frac{4\beta}{\beta'}},\label{ACOT3}
\end{equation}
for all $\varepsilon< \varepsilon_{0}=2^{-\frac{\beta'}{4}}$.

For any $i=0,\dots,n-1$, set
\begin{equation*}
F_{i}=\left\{\sup_{s_{i}\leq u \leq s_{i+1}}\left|
Y_{u}\left(u\right) \right| \leq \left(1+
\sqrt{\alpha_{1}}\right)\varepsilon^{r},\; \inf_{\left(u,\lambda
\right)\in \left[
s_{i},s_{i+1}\right]^{2}}\Upsilon_{u}\left(\lambda\right)\geq
\frac{\alpha_{2}}{2}\varepsilon ^{\frac{3}{2}}\right\}.
\end{equation*}
Using  $\left(\ref{ACOT3}\right)$, we prove that $D_{i} \subset
F_{i}$ , for all $\varepsilon <\varepsilon _{1} \wedge
\varepsilon_{0}$, for some $\varepsilon_{1}>0$.

Indeed,  the triangular inequality  and the estimate
$\left(\ref{ACOT3}\right)$ implies that, for any $\lambda\in[s_i,s_{i+1}]$, on $D_{i}$
$$
\Upsilon_{u}\left(\lambda\right) \geq
\Upsilon_{u}\left(u\right)-Cs^{\beta
}\varepsilon^{\frac{4\beta}{\beta'}} \geq
\alpha_{2}\varepsilon^{\frac{3}{2}}-Cs^{\beta
}\varepsilon^{\frac{4\beta}{\beta'}}.$$

Since $\beta'<\beta$, we have $\frac{4\beta }{\beta'}-\frac{3
}{2}>0$. Choosing $ \varepsilon_{1}=\left( \frac{\alpha_{2}}{2C
s^{\beta }}\right) ^{\frac{1}{^{\frac{4\beta
}{\beta'}-\frac{3}{2}}} }$, we obtain $
\Upsilon_{u}\left(\lambda\right) \geq
\frac{\alpha_{2}}{2}\varepsilon ^{\frac{3}{2}}$, for all
$\varepsilon <\varepsilon_{1} \wedge \varepsilon_{0}$ and for any
$(u,\lambda)\in \left[s_{i},s_{i+1}\right]^{2}$.

Thus, we have now reduced the proof to estimate
$\mathbf{P}\left\{ F_{i}\right\}$.

We shall only consider the case $i=0$. In fact, it will become clear
from the proof that the arguments depend only
on the length of the interval $\left[ s_{i},s_{i+1}\right]$.

We shall prove the existence of the
$\arg \sup$ associated to $Y_{\cdot}\left(\lambda\right)$.
This is done following the same arguments as in \cite{CM02}, Proposition 3.2,
that is, using Girsanov's theorem. Our setting needs a more general version
of this theorem than the one used in \cite{CM02}. For instance, we can apply Theorem 35 in \cite{Pr04}, p. 132.
With this, on a new probability space, the semimartingale $Y_{\cdot}\left(\lambda\right)$
is transformed into a local martingale and then, by a standard time change, into
a Brownian motion, for which the $\arg \sup$ does exist.
We only give an sketch of this procedure, since it is very similar as in \cite{CM02}.

Applying a Girsanov transformation needs to work on
the whole probability space $\Omega$, and not only on $F_{0}$.
For this reason, we have to modify the process
$Y_{\cdot}\left(\lambda\right)$, as follows.
Define
\begin{equation*}
\phi(z)=\left\{
\begin{array}{c}
z\;\;\ \ \ \ \ \ \ \ \ \text{if\ }\left| z\right|
\leq\mathcal{K}, \\
\pm 2\mathcal{K} \;\;\ \ \ \ \ \text{if\ }\left| z\right|
>3\mathcal{K},
\end{array}
\right. \;
\end{equation*}
otherwise, $\left|\phi\right|$ is bounded by $2\mathcal{K}$ and with
derivative bounded by 1;
\begin{equation*}
\psi_{\varepsilon }(z)=\left\{
\begin{array}{c}
1\;\;\ \ \ \ \ \ \ \ \ \ \text{if}\;\left| z\right| \geq
\frac{\alpha_{2}}{2}\varepsilon^{
\frac{3}{2}}, \\
\\
0\;\ \ \ \ \ \ \ \ \ \ \ \text{if}\;\left| z\right| <
\frac{\alpha_{2}}{4}\varepsilon^{ \frac{3}{2}},
\end{array}
\right. \;\;
\end{equation*}
$\psi_{\varepsilon}\;$is even and non-decreasing on
$[0,\infty)$.
Set
\begin{align}
\overline{\Phi}_{\cdot}\left(\ast\right)
&=\phi\left(\Phi_{\cdot}\left(\ast\right)\right), \label{coefbv}\\
\overline{\Psi}_{\cdot}^{j}\left(\ast\right)
&=\phi\left(\Psi_{\cdot}^{j}\left(\ast\right)\right)
\psi_{\varepsilon}\left(\Upsilon_{\cdot}\left(\ast\right)\right),\;1\leq
j\leq m, \label{coefsi1}\\
\overline{\Psi}_{\cdot}^{m+1}\left(\ast\right)
&=\sqrt{\frac{\alpha_{2}}{2}\varepsilon^{\frac{3}{2}}}
\left[1-\psi_{\varepsilon}\left(\Upsilon_{\cdot}\left(\ast\right)\right)\right].
\label{coefsi2}
\end{align}
Consider a Brownian motion $\widetilde{M}^{m+1}$, independent of
 $\left(\widetilde{M}^1,\ldots,\widetilde{M}^{m}\right)$.
Let
\begin{equation}
\overline{Y}_{u}(\lambda)=Y_{0}(\lambda)
+\sum_{j=1}^{m+1}\int_{0}^{u}{\overline{\Psi}}_{\eta}^{j}
\left(\lambda\right)d\widetilde{M}_{\eta}^{j}+\int_{0}^{u}\overline{\Phi}
_{\eta}\left(\lambda\right)d\eta \label{Solmod}.
\end{equation}
This is the analogue of (3.12) in \cite{CM02}.

Observe that if
$\omega \in F_{0}$ for  $u\leq s_{1}$, then
$\overline{Y}_{u}\left(\lambda,\omega
\right)=Y_{u}\left(\lambda,\omega\right)$.
 Hence, $\overline{Y}_{\cdot}\left(\lambda\right)$ is a
modification of $Y_{\cdot}\left(\lambda\right)$.

We define the modification of $M_{u}\left(\lambda\right)$ by
$
\overline{M}_{u}\left(\lambda\right):=\sum_{j=1}^{m+1}\int_{0}^{u}{\overline{\Psi}}_{\eta}^{j}
\left(\lambda\right)d\widetilde{M}_{\eta}^{j}$,
with quadratic variation
$
\left \langle \overline{M}_{\cdot}\left(\lambda\right) \right
\rangle_{u}=\int_{0}^{u}\overline{\Upsilon}_{\eta}\left(\lambda\right)d\eta$,
where
\begin{equation}
\overline{\Upsilon}_{\eta}(\lambda)=\sum_{j,k=1}^{m}\overline{\Psi}_{\eta}^{j}\left(\lambda\right)
\overline{\Psi}_{\eta}^{k}\left(\lambda\right)\Theta_{\eta}^{j,k}+
\left(\overline{\Psi}_{\eta}^{m+1}\right)^{2}\left(\lambda\right).\label{defPsimod}
\end{equation}
One can check that for all
$\eta$.
\begin{equation}
\label{lower}
\overline{\Upsilon}_{\eta}\left(\lambda\right)\geq
\frac{\alpha_{2}}{8}\varepsilon^{\frac{3}{2}}.
\end{equation}
This is the crucial fact ensuring the existence of a  probability measure
$\widetilde{\mathbf{P}}$, equivalent to ${\mathbf{P}}$, such that on each $\mathcal{F}_{s}$,
\begin{equation*}
\frac{d\widetilde{\mathbf{P}}}{d\mathbf{P}}
=\exp \left\{
-\int_{0}^{s}\frac{\overline{\Phi}_{\eta}\left(\lambda\right)}
{\overline{\Upsilon}_{\eta}\left(\lambda\right)}d\overline{M}_{\eta}\left(\lambda\right)
-\frac{1}{2}\int_{0}^{s}\frac{\overline{\Phi}_{\eta}^{2}\left(\lambda\right)}
{\overline{\Upsilon}^{2}_{\eta}\left(\lambda\right)}d\left \langle
\overline{M}_{\cdot}\left(\lambda\right)\right\rangle_\eta\right\},
\end{equation*}
and such that

\begin{equation*}
\overline{\mathbf{M}}_{u}\left(\lambda\right)=
\overline{M}_{u}\left(\lambda\right)+\int_{0}^{u}\overline{\Phi}_{\eta}\left(\lambda\right)d\eta,
\end{equation*}
is a $\widetilde{\mathbf{P}}$ local martingale, (see
$\cite{Pr04}$).
Thus, $\widetilde{\mathbf{P}}$-a.s.,
$
\overline{Y}_{u}(\lambda)=Y_{0}(\lambda)+
\overline{\mathbf{M}}_{u}\left(\lambda\right)$,
is a local martingale.

Set
\begin{align*}
\mathbf{A}_{u}(\lambda)&=\left \langle
\overline{Y}_{\cdot}(\lambda)\right \rangle_{u}
=\int_{0}^{u}\overline{\Upsilon}_{\eta}\left(\lambda\right)d\eta,\\
\mathbf{T}_{u}(\lambda)&=\inf \left\{ \eta\geq
0,\;\mathbf{A}_{\eta}(\lambda)\geq u\right\} .
\end{align*}
Then, there exists a $\widetilde{\mathbf{P}}$-Brownian motion $B$
such that for all $u$
$
\overline{Y}_{u}(\lambda)=Y_{0}(\lambda)+
B_{\mathbf{A}_{u}(\lambda)}$,
that is,
$\overline{Y}_{\mathbf{T}_{u}(\lambda)}(\lambda)=Y_{0}(\lambda)+
B_{u}$ $\widetilde{\mathbf{P}}$-a.s.,
(see \cite{KS91}, page 174, for details).

Let $S_{1}=\frac{\alpha_{2}\varepsilon
^{\frac{3}{2}}}{8}s_{1}\geq \frac{\alpha_{2}}{8}s\varepsilon
^{\frac{3}{2}+\frac{4}{\beta'}}$. Using
$(\ref{lower})$, we can prove
\begin{equation}
\label{1111}
u\frac{\alpha_{2}}{8}\varepsilon^{\frac{3}{2}}\le \mathbf{A}_{u}(\lambda)
\leq 2m^{2}\mathcal{K}^{3}u.
\end{equation}
Consequently,
$
\frac{u}{2m^{2}\mathcal{K}^{3}}\leq \mathbf{T}_{u}(\lambda)\leq
\frac{8u}{\alpha_{2}\varepsilon ^{\frac{3}{2}}}$.
Hence,
\begin{equation}
\frac{\alpha_{2}s\varepsilon
^{\frac{3}{2}+\frac{4}{\beta'}}}{16m^{2}\mathcal{K}^{3}}\leq
\frac{S_{1}}{2m^{2}\mathcal{K}^{3}}\leq
\mathbf{T}_{S_{1}}(\lambda)\leq
\frac{8S_{1}}{\alpha_{2}\varepsilon
^{\frac{3}{2}}}=s_{1}.\label{Vt}
\end{equation}

Set $\Pi_{1}(\lambda):=\mathbf{T}_{S_{1}}(\lambda)$.
The supremum of the absolute value of a linear
Brownian motion on a deterministic time interval is a.s. attained
at a single time. Thus,
for all $\lambda\in\left[ 0,s_{1}\right]$, $\widetilde{\mathbf{P}}$-a.s., there
exists an unique (random) time such that
\begin{equation}
\eta_{1}(\omega,\lambda)=\underset{0\leq \eta\leq
\Pi_{1}(\lambda)}{\arg \sup } \left|
\overline{Y}_{\eta}(\lambda)\right|. \label{argsup}
\end{equation}
Since $\mathbf{P}$ and $\widetilde{\mathbf{P}}$ are equivalent
 and $\Pi_{1}(\lambda)\leq s_{1}$,
$(\ref{argsup})$ holds $\mathbf{P} $-a.s. Thus, we have
$
\eta_{1}(\omega,\lambda)=\underset{0\leq \eta\leq
\Pi_{1}(\lambda)}{\arg \sup } \left|
\overline{Y}_{\eta}(\lambda)\right| \;\;\mathbf{P}\text{-a.s.}$, for all
$\lambda \in \left[ 0,s_{1}\right]$.

The final step consists of proving some control on the modified process
$\overline{Y}_{\cdot}(\lambda)$, because we only have the existence
of the $\arg\sup$ for this process.

Let $\theta_{1}=\frac{\alpha_{2}s\varepsilon
^{\frac{3}{2}+\frac{4}{\beta'}}}{16m^{2}\mathcal{K}^{3}}$.
Owing to (\ref{Vt}),
\begin{equation}
\sup_{0\leq \eta\leq \theta_{1}}\left|
\overline{Y}_{\eta}(\lambda)\right| \leq \sup_{0\leq \eta\leq
\Pi_{1}(\lambda)}\left| \overline{Y}_{\eta}(\lambda)\right|.
\label{desigualtat}
\end{equation}
The term $\sup_{0\leq \eta\leq \theta_{1}}\left|\overline{Y}_{\eta}(\lambda)\right|$
can be estimated  in a similar manner as in
$\cite{No86}$. We consider the It\^{o} formula for
$Y_{\theta_{1}}^{2}(\lambda)$,
\begin{equation}
\overline{Y}_{\theta_{1}}^{2}(\lambda
)=\overline{Y}_{0}^{2}(\lambda
)+2\int_{0}^{\theta_{1}}\overline{Y}_{\eta}(\lambda)d\overline{M}_{\eta}\left(\lambda\right)+2\int_{0}^{\theta_{1}}\overline{Y}_{\eta}
\left(\lambda\right)\overline{\Phi}_{\eta}\left(\lambda\right)d\eta+
\int_{0}^{\theta_{1}}\overline{\Upsilon}_{\eta}\left(\lambda\right)d\eta.\label{itoformulaqr}
\end{equation}

Fix $r^{\prime}
>\frac{2}{\beta^{\prime}}+\frac{11}{4}$. We have the following facts concerning
$\overline{Y}_{\cdot}(\lambda)$,
\smallskip

{\bf Fact 3} {\sl There exists $\varepsilon _{0}^{\prime}>0$ such
that for all $\varepsilon <\varepsilon _{0}^{\prime}$
\begin{equation*}
\mathbf{P}\left( \sup_{0\leq \eta\leq \theta_{1}}\left|
\overline{Y}_{\eta}(\lambda)\right| \leq
\left(2+\sqrt{\alpha_{1}}\right)\varepsilon ^{r^{\prime }}\right)
\leq 2\exp \left( -\frac{1}{2}\varepsilon^{-1}\right).
\end{equation*}
}

{\bf Fact 4} {\sl There exists $\varepsilon _{0}^{\prime
\prime}>0$ such that for all $\varepsilon <\varepsilon_{0}^{\prime
\prime}$
\begin{equation*}
\mathbf{P}\left( \exists \;\lambda \in \left[ 0,s_{1}\right]
,\;\sup_{0\leq \eta\leq \theta_{1}}\left|
\overline{Y}_{\eta}(\lambda)\right| \leq
\left(1+\sqrt{\alpha_{1}}\right)\varepsilon ^{r'}\right) \leq
c_{p}\varepsilon ^{p},
\end{equation*}
for all $p \geq 2$.
}

For their proofs, we can follow Lemma 3.21 and 3.22 of \cite{CM02},
respectively.

Let us return to $Y_{\cdot}\left(\lambda\right)$.
The processes $Y_u(\lambda)$, ${\bf A}_u(\lambda)$ and ${\bf T}_u(\lambda)$, have
jointly continuous sample paths in $(u,\lambda)$, a.s.
Thus, for $\omega \in F_{0}$, for all $\left(u,\lambda\right) \in \left[
0,s_{1}\right]^{2}$,
$\overline{Y}_{u}(\lambda,\omega
)=Y_{u}(\lambda,\omega)$.

Consider the set
\begin{equation*}
G_{0}=F_{0}\cap \left\{ \forall \;\lambda \in \left[
0,s_{1}\right] ,\;\sup_{0\leq \eta\leq
\Pi_{1}\left(\lambda\right)}\left| {Y}_{\eta}(\lambda)\right| >
\left(1+\sqrt{\alpha_{1}}\right)\varepsilon ^{r'}\right\}.
\end{equation*}
As in \cite{CM02}, one checks that $G_0=\emptyset$.
Consequently, with the result stated in Fact 4, we end the
proof of the lemma. \hfill\qed.

\section{Proof of Theorem \ref{t3.1}}
\label{s4}

To simplify the notation, we omit the summation sign on repeated indices.
Using the classical approach going back to Malliavin -and also by Bismut, Ikeda and
Watanabe, Stroock, Bouleau and Hirsch, etc.- (see for instance \cite{Nu95}), the proof of
the theorem consists of checking that
\begin{itemize}
\item[(i)] $X_z^i\in\mathbb{D}^{\infty }$, for all
$1\leq i\leq m$,
\item[(ii)]
$\det C_{z}^{-1}\in L^{p}$, for all $p\geq 2$, where $C_z$ is the $m\times m$ matrix whose entries are
\begin{equation*}
C_z^{ij}=
\int_{R_{z}}(D_{r}^{l}X_{z}^{i}) (D_{r}^{l}X_{z}^{j}) dr.
\end{equation*}
\end{itemize}

Proving (i) is strightforward. Indeed, due to the condition (h2), one can
proceed as in \cite{NS85}, Proposition 3.3.

To prove (ii), it suffices to show that for any $p\ge 2$ there exists $\varepsilon_0(p)$ such that for every $\varepsilon
\leq \varepsilon_{0}\left(p\right)$
\begin{equation}
\label{3.1}
\sup_{|v|=1}\textbf{P}\left\{v^{T}C_zv \leq \varepsilon\right\}\leq
\varepsilon^{p}
\end{equation}
(see e.g. \cite{Nu95})

We next give some preparations for the proof of (\ref{3.1}), valid in each one of the set of assumptions of Theorem \ref{t3.1}.
Following similar arguments as in \cite{NS85}, pp. 585-586, we write
$$v^TCv
=\sum_{l=1}^{d}\int_{R_{z}}\left
\langle v,\xi \left(r,z\right)A_{l}\left(r,X_{r}\right) \right
\rangle ^{2}dr,$$
where $v \in \R^{m}$,
and for  $r\in R_{S,T}$,
$\xi_{j}^{i}(r,z),\;r\preceq z$, $1 \leq i,j \leq m$,
is the solution to
\begin{equation}
\xi_{j}^{i}(r,z)=\delta_{j}^i+\int_{[r,z]}
\partial_{k}^{x}A_{l}^{i}\left(u,X_{u}\right)\xi_{j}^{k}(r,u)dW_{u}^{l}
+\int_{[r,z]}
\partial_{k}^{x}A_{0}^{i}\left(u,X_{u}\right)\xi_{j}^{k}(r,u)du.\label{firstvar}
\end{equation}
Then,
for any  $v\in S^{m-1}$, $0 < \varepsilon < 1$,  $0 < \mu <
1$, we have
$$
\mathbf{P}\left\{ v^{T}C_zv\leq \varepsilon \right\}\leq
\mathbf{P}(E_0)+\mathbf{P}(B),$$
where
\begin{align*}
E_0&=\left\{ \sum_{l=1}^{d}\int_{0}^{s}\left\langle
v,A_{l}\left( \eta ,t,X_{\eta ,t}\right) \right\rangle^{2}d\eta
\leq 4\varepsilon ^{\mu }\right\},\\
B&=\left\{\sum_{l=1}^{d}\int_{0}^{s}\int_{t-\varepsilon ^{1-\mu
}}^{t}\left\langle v,A_{l}\left( \eta,t,X_{\eta,t}\right) -\xi
\left( \eta, \tau,s,t\right) A_{l}\left(
\eta,\tau,X_{\eta,\tau}\right)\right\rangle ^{2}d\eta d\tau
>\varepsilon \right\}.
\end{align*}
Combining results of \cite{NS85} and the assumption (h1), we can obtain
an upper bound like in (\ref{3.1}) for $P(B)$, as follows.
Applying Chebyshev and  Cauchy-Schwarz inequalities yields
\begin{align*}
 &\textbf{P}\left\{B\right\}  \leq \varepsilon^{-\mu q}C\\
&\quad\times\sup_{\left(\eta,\tau\right)\in \left[0,s\right]\x
\left[t-\varepsilon^{1-\mu},t\right]}\max_{1 \leq l \leq
d}\mathbb{E}\left \| A_{l}\left( \eta,t,X_{\eta,t}\right)-\xi
\left( \eta, \tau,s,t\right) A_{l}\left(
\eta,\tau,X_{\eta,\tau}\right) \right \|^{2q}.
\end{align*}
Clearly,
\begin{align*}
& \mathbb{E}\left \| A_{l}\left( \eta,t,X_{\eta,t}\right)-\xi
\left( \eta, \tau,s,t\right) A_{l}\left(
\eta,\tau,X_{\eta,\tau}\right) \right \|^{2q}\\
& \leq C\mathbb{E}\left \| A_{l}\left(
\eta,t,X_{\eta,t}\right)-A_{l}\left(
\eta,\tau,X_{\eta,t}\right) \right \|^{2q} \\
& +C\mathbb{E}\left \| A_{l}\left(
\eta,\tau,X_{\eta,t}\right)-A_{l}\left(
\eta,\tau,X_{\eta,\tau}\right) \right \|^{2q}\\
& +C\mathbb{E}\left \| A_{l}\left(
\eta,\tau,X_{\eta,\tau}\right)-\xi \left( \eta, \tau,s,t\right)
A_{l}\left( \eta,\tau,X_{\eta,\tau}\right)  \right \|^{2q}.
\end{align*}
From (h1), we have
\begin{equation}
\mathbb{E}\left \| A_{l}\left(
\eta,t,X_{\eta,t}\right)-A_{l}\left( \eta,\tau,X_{\eta,t}\right)
\right \|^{2q} \leq K_{\gamma}^{2q}|t-\tau|^{2q\gamma},
\label{term1}
\end{equation}
By the mean value theorem in the spatial component and Lemma 3.1
in $\cite{NS85}$ applied to $(\ref{Itoeq})$,
\begin{equation}
\mathbb{E}\left \| A_{l}\left(
\eta,\tau,X_{\eta,t}\right)-A_{l}\left(
\eta,\tau,X_{\eta,\tau}\right) \right \|^{2q} \leq
C |t-\tau|^q, \label{term2}
\end{equation}
Finally, the Cauchy-Schwarz inequality and Lemma 3.1 in
$\cite{NS85}$ applied to $(\ref{firstvar})$ yield
\begin{align*}
& \mathbb{E}\left \| A_{l}\left(
\eta,\tau,X_{\eta,\tau}\right)-\xi \left( \eta, \tau,s,t\right)
A_{l}\left( \eta,\tau,X_{\eta,\tau}\right)  \right \|^{2q} \\
&=\mathbb{E}\left \| \left(I_{m}-\xi \left( \eta,
\tau,s,t\right)\right) A_{l}\left(
\eta,\tau,X_{\eta,\tau}\right)\right \|^{2q} \\
&\leq C \left|\left (s-\eta\right)\left
(t-\tau\right)\right|^q.
\end{align*}
Consequently,
\begin{align*}
\textbf{P}\left \{B\right\} \leq
C \varepsilon^{\left\{(1-2\mu)\wedge \left
(2\gamma-(2\gamma+1)\mu\right)\right\}q}.
\end{align*}
If $\mu < \frac{1}{2} \wedge \frac{2\gamma}{2\gamma+1}$ then
$(1-2\mu)\wedge \left (2\gamma-(2\gamma+1)\mu\right)>0$.
Thus, in this case, for all $p=\left\{(1-2\mu)\wedge \left
(2\gamma-(2\gamma+1)\mu\right)\right\}q \geq 2$,
\begin{equation*}
\mathbf{P}\left\{B\right\} \leq
C(K,K_{\gamma},S,T,d,p)\varepsilon^p.
\end{equation*}

We devote the remaining of the section to estimate the term
$\mathbf{P}(E_0)$.

Consider the stopping time $\mathcal{S}$ with respect to
the family of $\sigma$-algebras
$\left\{\mathcal{F}_{\eta,t},\;\eta \geq 0 \right\}$ (associated to the Brownian sheet) defined as
\begin{equation*}
\mathcal{S}=\inf \left\{ \eta \geq 0:\; \sup_{\theta \leq \eta;\;
\tau\le t} \left| X_{\theta,\tau}-x_{0}\right|\geq R\right\}
\wedge s \wedge s_{0}.
\end{equation*}
For any $k=0,1,\ldots,N$, $v \in S^{m-1}$, define
\begin{equation*}
E_k:=\left \{\sum_{V \in { \Sigma }_{k}}\int_{0}^{s} \left \langle
v,V\left(\eta,t,X_{\eta,t}\right)\right \rangle^2d\eta \leq
4\varepsilon^{m(k)} \right\}
\end{equation*}
and $E=\overset{N}{\underset{k=0}{\cap }}E_{k}$,
where the $m(k)$ are positive constants to be fixed later.

As usually, for $N\ge 1$ we consider the decomposition\\
$E_{0}\subset \left( E_{0}\cap E_{1}^{c}\right) \cup \left(
E_{1}\cap E_{2}^{c}\right) \cup \cdots \cup \left( E_{N-1}\cap
E_{N}^{c}\right) \cup E$,
yielding
\begin{equation}
\label{decomposition}
\mathbf{P}(E_{0})\leq \mathbf{P}(E)+\sum_{k=0}^{N-1}\mathbf{P}\left( E_{k}\cap
E_{k+1}^{c}\right).
\end{equation}
We are going to estimate each term of this sum under the different set of assumptions
of the theorem.

Assume (h3); we  shall prove that
$\mathbf{P}(E)\leq C \varepsilon ^{p}$,
for any $\varepsilon>0$ small enough, where $C$ depends on $t$, $p$, $m$, $K$ and $R$.

We write
$\mathbf{P}\left( E\right) \leq \mathbf{P}\left( E\cap \left\{
\mathcal{S}\geq \varepsilon^{\varsigma }\right\} \right)
+\mathbf{P}\left( \mathcal{S}<\varepsilon ^{\varsigma }\right),
$
with $0<\varsigma <M=\underset{k=0,\ldots ,N}{\min }m(k)$.

For $\varepsilon$ small enough, the set
 $E\cap \left\{ \mathcal{S}\geq \varepsilon^{\varsigma}\right\}$ is empty.
 In fact, assume $\mathcal{S}\geq
\varepsilon^{\varsigma}$; from $\left(\ref{Hör''}\right)$ we have
\begin{equation*}
\sum_{k=0}^{N}\sum_{V\in { \Sigma }_{k}}\int_{0}^{s}\left\langle
v,V\left(\eta,t,X_{\eta,t}\right) \right\rangle ^{2}d\eta \geq
c_N \varepsilon ^{\varsigma }.
\end{equation*}

On the set $E$,
\begin{equation*}
\sum_{k=0}^{N}\sum_{V\in { \Sigma }_{k}}\int_{0}^{s}\left\langle
v,V\left(\eta,t,X_{\eta,t}\right) \right\rangle ^{2}d\eta \leq
4\left( N+1\right) \varepsilon ^{M}.
\end{equation*}

Thus, for all $\varepsilon <\varepsilon_{0}=\left(
\frac{c_{N}}{4\left(N+1\right)}\right)^{\frac{1}{M-\varsigma}}$,
 $E\cap \left\{ \mathcal{S}\geq
\varepsilon^{\varsigma}\right\} =\emptyset$.

Applying  Chebychev, Burkholder-Davis-Gundy and
H\"{o}lder inequalities, Gronwall's lemma and Lemma 3.1 in
$\cite{NS85}$ to $(\ref{Itoeq})$, yields
\begin{align*}
\mathbf{P}\left( \mathcal{S}<\varepsilon ^{\varsigma }\right)\leq
C \varepsilon ^{\frac{p^{\prime }\varsigma
}{2}},
\end{align*}
for all $p^{\prime }\geq 2.$
Thus, taking $p=\frac{p^{\prime }\varsigma}{2}\geq 2,$ for all
$\varepsilon < \varepsilon_0$
\begin{equation*}
\mathbf{P}\left( E\right) \leq \mathbf{P}\left(
\mathcal{S}<\varepsilon ^{\varsigma }\right) \leq C \varepsilon
^{p},
\end{equation*}
where $C$ depends on $p$, $t$, $m$, $K$ and $R$.

Assume $N=0$. Then $E=E_0$;    
therefore, for the {\bf elliptic case} the proof is complete.
\medskip

{\bf The smooth case}

In addition to  the results proved so far, we have to
study  $\mathbf{P}(E_{k}\cap E_{k+1}^{c})$,
$k=0,\ldots,N-1$. Clearly, $E_{k}\cap E_{k+1}^{c}\subset \cup_{V\in \Sigma_k} B$, where
\begin{multline*}
B=\left \{\int_{0}^{s}\left\langle
v,V\left(\eta,t,X_{\eta,t}\right)
\right\rangle ^{2}d\eta\leq 4\varepsilon ^{m(k)}, \right. \\
\left.\sum_{l=1}^{d}\int_{0}^{s}\left\langle
v,A_{l}^{\nabla}V\left( \eta,t,X_{\eta,t}\right) \right\rangle
^{2}d\eta>\frac{4\varepsilon^{m(k+1)}}{\Gamma}\right \},
\end{multline*}
and $\Gamma=\underset{k=0,\ldots,N-1}\max
card\left(\Sigma_{k}\right)$.

Moreover, $B\subset B^1\cup B^2$, with
\begin{multline*}
B^{1}:=\left \{ \int_{0}^{s} \left \langle
v,V\left(\eta,t,X_{\eta,t}\right) \right \rangle^2 d\eta \leq
4\varepsilon^{m(k)}, \right. \\
\left.\sum_{l=1}^{d}\int_{0}^{s}\int_{0}^{t}\left \langle
v,A_{l}\left(\eta,\tau,X_{\eta,\tau}\right)^{\nabla}
V\left(\eta,t,X_{\eta,t}\right)\right\rangle^2 d\eta d\tau \geq
\frac{\varepsilon^{\vartheta m(k+1)}}{\Gamma} \right \},
\end{multline*}
\begin{multline*}
B^{2}:=\left \{ \sum_{l=1}^{d}\int_{0}^{s}\int_{0}^{t}\left
\langle v,A_{l}\left(\eta,\tau,X_{\eta,\tau}\right)^{\nabla}
V\left(\eta,t,X_{\eta,t}\right)\right\rangle^2 d\eta d\tau <
\frac{\varepsilon^{\vartheta m(k+1)}}{\Gamma} , \right. \\
\left.\sum_{l=1}^{d}\int_{0}^{s}\left\langle
v,A_{l}^{\nabla}V\left( \eta,t,X_{\eta,t}\right) \right\rangle
^{2}d\eta > \frac{4\varepsilon^{m(k+1)}}{\Gamma} \right \},
\end{multline*}
where $\vartheta$ is an arbitrary positive real constant.

Let us estimate $\mathbf{P}(B^{2})$. Firstly, one
easily checks that $B^2\subset {B^{2}}^{\prime}$ with
\begin{multline*}
{B^{2}}^{\prime}= \left \{
\sum_{l=1}^{d}\int_{0}^{s}\int_{t-\varepsilon^{\left(\vartheta-1\right)m(k+1)}}^{t}\left\langle
v,\left\{A_{l}\left(
\eta,t,X_{\eta,t}\right)-A_{l}\left(\eta,\tau,X_{\eta,\tau}\right)\right\}^\nabla\right.\right. \\
\left.\left.  V\left( \eta,t,X_{\eta,t}\right)
\right\rangle^2d\eta d\tau > \frac{\varepsilon^{\vartheta
m(k+1)}}{\Gamma}\right\}.
\end{multline*}

Applying Chebyshev's and the Cauchy-Schwarz inequalities, yields
\begin{align*}
&\mathbf{P}({B^{2}}^{\prime}) \leq
C \varepsilon^{-qm(k+1)}\\
&\quad\times\underset{\tau
\in
\left[t-\varepsilon^{\left(\vartheta-1\right)m(k+1)},t\right]}{\sup_{\eta
\in \left[0,s\right]}}\underset{1 \leq j \leq m}{\max_{1 \leq l
\leq d}}\Big \{\mathbb{E} \left \|A_{l}\left(
\eta,t,X_{\eta,t}\right)-A_{l}\left(\eta,\tau,X_{\eta,\tau}\right)\right\|^{4q}\\
&\quad \times \mathbb{E}\left\|\partial_{j}^{x}V\left(
\eta,t,X_{\eta,t}\right)\right\|^{4q} \Big\}^{\frac{1}{2}}
\end{align*}
We have
$
\mathbb{E} \left
\|\partial_{j}^{x}V\left(\eta,t,X_{\eta,t}\right)\right\|^{4q}
\leq K^{4q}
$. Moreover, owing to
 $(\ref{term1})$ and $(\ref{term2})$,
\begin{align*}
&\mathbb{E}\left\|A_{l}\left(\eta,t,X_{\eta,t}
\right)-A_{l}\left(\eta,\tau,X_{\eta,\tau}
\right)\right\|^{4q} \\
&\leq 2^{4q-1}\left
(K_{\gamma}^{4q}|t-\tau|^{4q\gamma}+C\left(K,S,T,q\right)|t-\tau|^{2q}\right).
\end{align*}
Thus,
$
\mathbf{P}(B^{2}) \leq
C\varepsilon^{\left\{\left(2\left(\vartheta-1\right)\gamma-1\right)\wedge
\left(\vartheta-2\right)\right\}qm(k+1)}.
$
Taking $\vartheta>\left(\frac{1}{2\gamma}+1\right) \vee 2$, we
obtain $\left\{\left(2\left(\vartheta -1\right)\gamma -1\right)
\wedge \left (\vartheta-2\right)\right\}>0$.
Hence, $\mathbf{P}(B^{2})\leq
C \varepsilon^{p}$, for all\\
$p=\left\{\left(2\left(\vartheta-1\right)\gamma-1\right)\wedge
\left(\vartheta-2\right)\right\}qm(k+1) \geq 2$.
\smallskip


It remains to study the term $\mathbf{P}(B^{1})$.
For fixed $t$, consider the one-parameter semimartingale $(X_{s,t}, s\in[0,S])$
with respect to the filtration $\left\{\mathcal{F}_{\eta,t},\;\eta \geq 0 \right\}$,

An application of the It\^{o} formula yields
\begin{align}
& \left \langle v,V\left(s,t,X_{s,t}\right)\right\rangle=\left
\langle v,V\left(0,t,x_{0}\right)\right\rangle+
\int_{0}^{s}\left \langle v,\partial^{\eta}V\left(\eta,t,X_{\eta,t}\right)\right\rangle d\eta \notag\\
&+\sum_{l=1}^{d}\int_{0}^{s}\int_{0}^{t}\left \langle
v,A_{l}\left(\eta,\tau,X_{\eta,\tau}\right)^{\nabla}
V\left(\eta,t,X_{\eta,t}\right)\right\rangle
dW_{\eta,\tau}^{l} \notag\\
&+\int_{0}^{s}\int_{0}^{t}\left \langle
v,A_{0}\left(\eta,\tau,X_{\eta,\tau}\right)^{\nabla}V\left(\eta,t,X_{\eta,t}\right)\right\rangle
d\eta d\tau \notag\\
&+\frac{1}{2}\int_{0}^{s}\int_{0}^{t}
\left\langle
v,\partial_{j}^{x}\partial_{k}^{x}V\left(\eta,t,X_{\eta,t}\right)A_{l}^{j}\left(\eta,\tau,X_{\eta,\tau}\right)A_{l}^{k}
\left(\eta,\tau,X_{\eta,\tau}\right)\right\rangle d\eta d\tau
\label{ItoformulaC1}.
\end{align}

We can now use the arguments of \cite{NS85}. More precisely, we apply  Lemma $4.2$ in \cite{NS85}
to the continuous semimartingale $Y_s:=\langle v,V(s,t,X_{s,t})\rangle$, which decomposition follows
from (\ref{ItoformulaC1}). This finishes the proof in the  regular case.
\smallskip

The estimate of $P(E_0)$ in the {\bf factorable case} can be obtained
using the same method as for the regular case. We skip the details
of the proof to avoid repetitions.
\medskip

{\bf The regular H\"older case}

Under this set of assumptions, the expression  (\ref{ItoformulaC1}) does not make sense,
because the vector fields $V(\eta,t,x)$ are not differentiable with respect to the variable $\eta$.
Instead, we consider the process $\left(\langle v, V(\lambda, t, X_{\eta,t}), \eta\in[0,s]\right)$,
where $V\in \Sigma_k$, $k=0,\dots, N-1$, and $\lambda \ge 0$, $t\in[0,T]$, are fixed, and apply
the It\^o formula. We obtain
\begin{align}
& \left \langle
v,V\left(\lambda,t,X_{s,t}\right)\right\rangle=\left \langle
v,V\left(\lambda,t,x_{0}\right)\right\rangle\notag\\
&+\sum_{l=1}^{d}\int_{0}^{s}\int_{0}^{t}\left \langle
v,A_{l}\left(\eta,\tau,X_{\eta,\tau}\right)^{\nabla}
V\left(\lambda,t,X_{\eta,t}\right)\right\rangle
dW_{\eta,\tau}^{l} \notag\\
&+\int_{0}^{s}\int_{0}^{t}\left \langle
v,A_{0}\left(\eta,\tau,X_{\eta,\tau}\right)^{\nabla}V\left(\lambda,t,X_{\eta,t}\right)\right\rangle
d\eta d\tau \notag\\
&+\frac{1}{2}\int_{0}^{s}\int_{0}^{t}\left \langle
v,\partial_{j}^{x}\partial_{k}^{x}V\left(\lambda,t,X_{\eta,t}\right)A_{l}^{j}\left(\eta,\tau,X_{\eta,\tau}\right)A_{l}^{k}
\left(\eta,\tau,X_{\eta,\tau}\right)\right\rangle d\eta
d\tau.\label{ItoformulaHC}
\end{align}
Set
\begin{align*}
Y_{s}\left(\lambda\right)&=\left \langle
v,V\left(\lambda,t,X_{s,t}\right)\right\rangle, \quad
Y_{0}\left(\lambda\right)=\left \langle v,V\left(\lambda,t,x_{0}\right)\right\rangle,\\
\Psi_{\eta}^{j} \left(\lambda\right)&=\left\langle v,
\partial_{j}^{x}V\left(\lambda,t,X_{\eta,t}\right)\right\rangle, \quad
\widetilde{M}^{j}_{s}=\int_{0}^{s}\int_{0}^{t}A_{l}^{j}\left(\eta,\tau,X_{\eta,\tau}\right)dW^{l}_{\eta,\tau},\\
\Theta_{\eta}^{j,k}&=\int_{0}^{t}\sum_{l=1}^{d}A_{l}^{j}\left(\eta,\tau,X_{\eta,\tau}\right)
A_{l}^{k}\left(\eta,\tau,X_{\eta,\tau}\right)d\tau,\\
\Phi_{\eta}\left(\lambda\right)&=\int_{0}^{t}\left \langle
v,A_{0}\left(\eta,\tau,X_{\eta,\tau}\right)^{\nabla}V\left(\lambda,t,X_{\eta,t}\right)\right\rangle
 d\tau \\
&+\frac{1}{2}\int_{0}^{t}
\left\langle
v,\partial_{j}^{x}\partial_{k}^{x}V\left(\lambda,t,X_{\eta,t}\right)A_{l}^{j}\left(\eta,\tau,X_{\eta,\tau}\right)A_{l}^{k}
\left(\eta,\tau,X_{\eta,\tau}\right)\right\rangle d\tau,
\end{align*}
and
$\alpha_{1}=4$, $\alpha_{2}=\frac{1}{\Gamma}$.

These processes satisfy the assumptions of Lemma \ref{l3.1}. Fix $\nu>\frac{3}{2\beta-1}$, then $3+2\nu<\rho$. Set $m(0)= \mu$ and $m(k)= \frac{\mu}{(\theta \rho)^k}$,
$k=0,\ldots,N$. With Lemma \ref{l3.1}, we obtain
the desired estimate for $P(B^1)$ and we finish the proof.

\medskip

{\bf The irregular H\"{o}lder case}

We apply Lemma \ref{l3.2} to the same processes as we did in the
regular H\"{o}lder case. Notice that the assumptions of this lemma
are satisfied. Fix
$\rho>\left(\frac{11}{2}+\frac{4}{\beta'}\right)\left(1+\frac{1}{\beta'}\right)$.
Set $m(0)= \mu$ and $m(k)= \frac{\mu}{(\vartheta \rho)^k}$,
$k=0,\ldots,N$. Then, for $\varepsilon$ small enough we obtain the
suitable estimate for $P(B^1)$, and therefore the proof is
complete. \hfill\qed
\bigskip

\noindent {\textbf{Acknowledgements}} The authors would like to
thank P. Cattiaux and L. Mesnager for crucial discussions on some
aspects of the proof of  Lemma \ref{l3.1}.

A part of this work has been done while the first author was
visiting the {\sl Centro di Recerca Matematica Ennio De Giorgi} at
Pisa. She would like to express her gratitude for the warm
hospitality and the financial support provided by the Centre.

\end{document}